\documentclass[12pt]{article}
\usepackage{latexsym}
\usepackage{amsfonts}
\def\R{\mathbf{R}}
\def\C{\mathbf{C}}
\def\bC{{\mathbf{\overline{C}}}}
\def\const{\mathrm{const}}
\def\Fr{\mathrm{Fr}\ }
\title{On the ``pits effect'' of Littlewood and Offord}
\author{Alexandre Eremenko\thanks{Supported by the NSF grants
DMS-0555279 and DMS-0244547.}$\;$ and Iossif Ostrovskii}
\date{October 13, 2007}
\begin{document}
\maketitle
\begin{abstract}
Asymptotic behavior of the entire functions
$$f(z)=\sum_{n=0}^\infty e^{2\pi i\alpha_n}z^n/n!,\quad\mbox{with real}\quad\alpha_n$$
is studied. It turns out that the Phragm\'en--Lindel\"of indicator
of such function is always non-negative, unless $f(z)=e^{az}$.
For special choice $\alpha_n=\alpha n^2$ with irrational $\alpha,$
the indicator is constant and $f$ has completely regular growth
in the sense of Levin--Pfluger.
Similar functions
of arbitrary order are also considered.

MSC Primary: 30D10, 30D15, 30B10.
\end{abstract}

In \cite{Nassif} Nassif studied (on Littlewood's suggestion)
the asymptotic behavior and the distribution of zeros
of the entire function
\begin{equation}
\label{11}
\sum_{n=0}^\infty e^{2\pi in^2\alpha}z^n/n!,
\end{equation}
with $\alpha=\sqrt{2}$.
This was continued by Littlewood \cite{Lit1,Lit2}, 
who considered generalizations to
Taylor series whose coefficients have smoothly varying
moduli and arguments of the form $\exp(2\pi i\alpha n^2)$,
where $\alpha$ is a quadratic irrationality.

Such functions behave similarly to random entire functions
previously studied by Levy \cite{Levy} and Littlewood and Offord \cite{LO},
in particular they display
the ``pits effect'' which Littlewood described as follows:

{\em ``If we erect an ordinate $|f(z)|$ at the point $z$ of
the $z$-plane, then the resulting
surface is an exponentially rapidly rising bowl, approximately
of revolution,
with exponentially small pits going down to the bottom.
The zeros of $f$, more generally the $w$-points where $f=w$,
all lie in the pits for $|z|>R(w)$. Finally
the pits are very uniformly distributed in direction, and as uniformly
distributed in distance as is compatible with the order $\rho$''}.

The earliest study of functions (\ref{11}) known to the authors
is the thesis of {\AA}lander \cite{AA} who considered the case
of rational $\alpha$. Levy \cite{Levy} used the results of Hardy and 
Littlewood on Diophantine approximation to prove the following. Let
$$M(r,f)=\max_{|z|=r}|f(z)|\quad\mbox{and}\quad
m_2^2(r,f)=\frac{1}{2\pi}\int_{-\pi}^{\pi}|f(re^{i\theta})|^2d\theta.$$
Then 
\begin{equation}
\label{ratio}
M(r,f)/m_2(r,f) \quad\mbox{is bounded}
\end{equation}
for $f$ of the form (\ref{11}), and
$\alpha$ satisfying a Diophantine condition. This is even stronger
regularity than random arguments of coefficients yield \cite{Levy,LO}.
Some other works where the function (\ref{11}) with various $\alpha$
was studied or used are 
\cite{E,E1,Mcintyre,Valiron}. 

Function (\ref{11}) is the unique analytic solution
of the functional equation 
\begin{equation}
\label{fae}
f'(z)=qf(q^2z),\quad\mbox{where}\quad q=e^{2\pi i\alpha},\quad\mbox{and}\quad
f(0)=1,
\end{equation}
which is a special case of the so-called ``pantograph equation''.
There is a large literature on this equation with real $q$,
see, for example, \cite{Langley,Iserles} and references there.

Recently
there was a renewed interest to the functions of the type (\ref{11})
because they arise as the limits as $q\to e^{2\pi i\alpha}$ of
the function of two variables
$$\sum_{n=0}^\infty q^{n^2}z^n/n!$$
which plays an important role in graph theory \cite{Tutte}
and statistical
mechanics \cite{Sokal}. This function is
the unique solution of (\ref{fae}), for all $q$ in the closed unit
disc.

In the present paper, we first
study arbitrary entire functions of the form
\begin{equation}
\label{1}
f(z)=\sum_{n=0}^\infty a_nz^n/n!,\quad\mbox{where}\quad |a_n|=1.
\end{equation}
Our Theorem~1 says that such functions
cannot decrease exponentially on
any ray, unless $f$ is an exponential.
This can be compared with a result of Rubel and Stolarski \cite{RS}
that there exist exactly five series of the form (\ref{1})
with $a_0=0,\;a_n=\pm1$
which are bounded on the negative ray.
Our second result, Theorem~2 shows that one cannot replace the
condition of exponential decrease in  Theorem~1 by boundedness
on a ray: there are infinitely many functions of the form
(\ref{1}) which tend to zero as $z\to\infty$
in the closed right half-plane.

In the second part of the paper, we consider the case
$\arg a_n=2\pi in^2\alpha$ with {\em any irrational}
$\alpha$. Theorem~3 shows that 
the qualitative picture of $|f(z)|$ is the same as described
by Littlewood, except that our estimate of
the size of the pits is worse than exponential.
In particular, we show that
$$\log|f(z)|=|z|+o(|z|),$$
outside some exceptional set of $z$.
According to the Levin--Pfluger theory \cite{L},
this behavior of $|f|$ has the following consequences
about the zeros $z_k$ of $f$:
\vspace{.1in}

The number $n(r,\theta_1,\theta_2)$ of zeros (counting
multiplicity) in the sector
$$\{ z: \theta_1<\arg z<\theta_2,\; |z|<r\}$$ satisfies
\begin{equation}
\label{p1}
n(r,\theta_1,\theta_2)=
\frac{\theta_2-\theta_1}{2\pi}(r+o(r))\quad\mbox{as}\quad
r\to\infty.
\end{equation}
Moreover, the
limit
\begin{equation}
\label{p2}
\lim_{R\to\infty}\sum_{|z_k|\leq R}\frac{1}{z_k}\quad\mbox{exists},
\end{equation}
where $z_k$ is the sequence of zeros of $f$. It is easy to see from
the Taylor series of $f$ that this limit equals $-q$.
\vspace{.1in}

Thus the Diophantine conditions used in
\cite{Levy,Valiron,Nassif} are unnecessary for the
qualitative picture of behavior of $|f|$, but with arbitrary
irrational $\alpha$ the results are less precise then those
where $\alpha$ satisfies a Diophantine condition.
Theorem~4 shows that Levy's result (\ref{ratio}) cannot be extended
to arbitrary irrational $\alpha$. 
Finally we prove 
a result similar to Theorem~3 where the condition $|a_n|=1$ is replaced by
a more flexible condition on the moduli of the coefficients
allowing the function to have any
order of growth.

We denote by
$$F(z)=\sum_{n=1}^\infty a_{n-1}z^{-n},$$
the Borel transform of $f$ in (\ref{1}) (terminology of \cite{L}).
Then $F$ has an analytic continuation from a neighborhood of
infinity to the region
$\bC\backslash K$, where $K$
is a convex compact set in the plane, which is called the conjugate
indicator
diagram. 
The indicator 
$$h_f(\theta):=\limsup_{r\to\infty}r^{-1}\log|f(re^{i\theta})|,\quad
|\theta|\leq\pi,$$
is the support function of the convex set
symmetric to $K$ with respect to the real axis.

We also consider the function
$$G(z)=\sum_{n=1}^\infty a_{n-1}z^n$$
analytic in the unit disc. Transition from $F$ to $G$ is by
the change of the variable $1/z$.
\vspace{.1in}

\noindent
{\bf P\'olya's theorem} (\cite[Appendix I, \S5]{L}).
{\em Suppose that $G$ has an analytic continuation from the
unit disc to infinity through some angle $|\arg z-\pi|<\delta.$
Then the coefficients
$a_n$ can be interpolated by an entire function $g$ of exponential
type such that the indicator diagram
of $g$ is contained in the horizontal strip
$|\Im z|\leq \pi-\delta.$ That is
$g(n)=a_n$ for $n\geq 1$, 
and $h_g(\pm\pi/2)\leq\pi-\delta.$}
\vspace{.1in}

\vspace{.1in}

\noindent
{\bf Carlson's theorem} (\cite[Ch. IV,
Intro.]{L}). {\em Suppose that the indicator diagram
of an entire function $g$ has width less than $2\pi$ in the direction
of the imaginary axis, that is $h_g(\pi/2)+h_g(-\pi/2)<2\pi$.
Then $g$  cannot vanish on the  positive integers,
unless $g=0$.}
\vspace{.1in}

\noindent
{\bf Theorem 1.} {\em Every entire function $f$ of the form $(\ref{1})$
has non-negative indicator,
unless $a_n=\const\cdot a^n$ for some $a$ on
the unit circle, in which case $f(z)=e^{az}.$}
\vspace{.1in}

By Borel's transform, this is equivalent to 
\vspace{.1in}

\noindent 
{\bf Theorem ${\mathbf{1'}}$.} {\em Let $G$ be as above. 
Then $G$ cannot have an analytic continuation to infinity through any
half-plane containing $0$, unless $a_n=\const\cdot a^n$ for some~$a$.}
\vspace{.1in} 

These two theorems give characterizations of the exponential
function and the geometric series, respectively, showing that their
behavior is quite exceptional. Another somewhat similar characterization
follows from the result in \cite{RS} mentioned above.
As a corollary from Theorem~1 we obtain
the following result of Carlson \cite{Carl}: If $z_n$ is the sequence
of zeros of $f$ as in (\ref{1}), then 
\begin{equation}
\label{carlson}
\sum_n1/|z_n|=\infty,
\end{equation}
unless $f$ is an exponential.
\vspace{.1in}

{\em Proof of Theorem $1'$.}
Suppose that $G$ has such an analytic continuation.
Replacing $z$ by $az$ with $|a|=1$ we achieve that $G$
has an analytic continuation to infinity through some left half-plane
of the form $\Re z <\epsilon$, where $\epsilon > 0$.

P\'olya's theorem then implies that $a_n=g(n)$ for some entire
function whose indicator diagram is contained in
the strip $|\Im z|<\pi/2-\delta$, for some $\delta >0$.
Consider the functions $$g_R(z)=(g(z)+\overline{g(\overline{z})})/2\quad
\mbox{and}
\quad g_I(z)=(g(z)-\overline{g(\overline{z})})/(2i).$$
On the real axis we have $g_R(x)=\Re g(x)$ and $g_I(x)=\Im g(x)$.
Consider the entire function
$$H=g_I^2+g_R^2.$$
Then at positive integers we have
$$H(n)=g_I^2(n)+g_R^2(n)=(\Re g(n))^2+(\Im g(n))^2=|a_n|^2=1.$$
So the function $H-1$ has zeros at all positive integers.
Its indicator diagram is contained in the strip
$$|\Im z|\leq \pi-2\delta,$$
(Squaring stretches the indicator diagram by a factor of 2,
and
the indicator diagram of the sum of
two functions is contained in the convex hull of the union of their diagrams).
Now, by Carlson's theorem, $H\equiv 1$, so 
\begin{equation}\label{pifagor}
g_I^2+g_R^2=1.
\end{equation} 
The general solution of this functional equation
in the class of entire functions
is $g_I=\cos\circ \phi,\; g_R=\sin\circ\phi$,
where $\phi$ is an entire function.
It is well-known and easy to see that for $g_I$ and $g_R$ to be of exponential
type, it is necessary and sufficient that $\phi(z)=cz+b$.
As $g_I$ and $g_R$ are real on the real line, we conclude that $c$ and $b$ are
real. Thus $a_n=\cos(cn+b)+i\sin(cn+b)=\const\cdot e^{icn}=\const\cdot a^n,$
as advertised.

An alternative way to derive the conclusion from (\ref{pifagor})
suggested by Katsnelson is
to notice that (\ref{pifagor}) implies 
\begin{equation}
\label{dwa}
|g(x)|\equiv 1\quad\mbox{for real}\quad x.
\end{equation}
The Symmetry Principle then implies that $g$ has no zeros (if $z_0$ is
a zero then $\overline{z}_0$ would be a pole). So $g$ is a function of
exponential type without zeros, so $g= \exp(icz)$,
where $c$ should be real by (\ref{dwa}).
\hfill$\Box$
\vspace{.1in}

As we already noticed, Theorem~1 implies (\ref{carlson}).
However it does not imply that the sequence of zeros has positive
density:
there exist functions of exponential type, even with constant indicator,
whose zeros have zero density.\footnote{
Valiron \cite[p. 415]{Valiron} erroneously asserted the contrary:
that for functions with constant indicator, zero cannot be
a Borel exceptional value.} To construct such examples,
take zeros of the form 
$$z_k=\left(e^{i\log\log(k+1)}-e^{i\log\log k}\right)^{-1},$$
and construct the canonical product $W$ of
genus one with such zeros.
It is not hard to show that 
the asymptotic behavior of this product will be
$$\log|W(re^{i\theta})|=(cr+o(r))\cos(\theta-\log\log r),\quad r\to\infty$$
outside of some small exceptional set, 
so the indicator $h_W$ is constant, while the density of
zeros is zero.
\vspace{.1in}

There exist entire functions of the form (\ref{1}), other than the
exponential, which are bounded 
in the left half-plane.
The simplest example is Hardy's generalization of $e^z$ defined
by the power series
$$E_{s,a}=\sum_{n=1}^\infty\frac{(n+a)^sz^n}{n!},\quad s\in\C,\quad
a>0.$$
For pure imaginary $s$, this series is of the form (\ref{1}).
Hardy \cite{Hardy} proved the asymptotic formula
$$E_{s,a}(z)
=z^se^z(1+o(1))+\frac{\Gamma(a)}{\Gamma(-s)(-z)^a\log(-z)}
(1+o(1)),$$
as $z\to\infty,\;|\arg z\pm\pi/2|<\epsilon,$
for every $\epsilon\in(0,\pi/2)$. This formula implies that
the functions $E_{s,a}$ with pure imaginary $s$
are bounded in the closed left half-plane.
For further results on Hardy's function, see \cite{O}.
\vspace{.1in}

\noindent
{\bf Theorem 2.} {\em Let $\psi$ be a real entire function
with the property
$$\psi(\zeta)=o(|\zeta|),\quad \zeta\to\infty$$
in every half-plane $\Re\zeta>c,\; c\in\R$. 
Then the function
$$f(z)=\sum_{n=0}^\infty\frac{e^{i\psi(n)}}{n!}z^n$$
is of the form $(\ref{1})$ and
for every $A>0$ and every $\epsilon>0$
we have 
\begin{equation}
\label{01}
|f(re^{i\phi})|=O(r^{-A}),\quad r\to\infty,
\end{equation}
uniformly for $|\phi-\pi|\leq\pi/2-\epsilon.$
}
\vspace{.1in}

{\em Proof.} We have the following integral representation:
\begin{equation}
\label{o1}
f(-z)=-\frac{1}{2\pi i}\int_{-A-i\infty}^{-A+i\infty}
\frac{\pi e^{i\psi(\zeta)}z^{\zeta}}{\Gamma(\zeta+1)\sin\pi\zeta}d\zeta,
=\frac{1}{2\pi i}\int_{-A-i\infty}^{-A+i\infty}e^{i\psi(\zeta)}z^{\zeta}
\Gamma(-\zeta)d\zeta,
\end{equation}
where $A>0$ is any positive number.
To obtain this representation, we notice that
that by Stirling's
formula, the modulus of the integrand does not exceed
$$|z|^{-\Re\zeta}\exp\left((-\pi/2+\phi+o(1))|\Im\zeta|\right),$$
as $|\zeta|\to \infty$ in every half-plane
of the form $\Re\zeta\geq -A$. Here $o(1)$ is independent of $z$.
Applying the residue formula to the rectangle
$$\{\zeta:-c<\Re\zeta<N+1/2,\, |\Im\zeta|<N+1/2\},$$
and letting $N$ tend to infinity, we obtain (\ref{o1}).
Now the same estimate of the integrand shows that (\ref{01})
holds. \hfill$\Box$

\vspace{.1in}

Theorem~1 implies that the indicator diagram
of a function of the form (\ref{1}),
other than an exponential, contains zero.
Theorem~2 shows that the indicator diagram of such a function
can be contained in a closed
half-plane.
It seems interesting to describe all possible indicator diagrams
that can occur for functions of the form (\ref{1}).
We have the following partial result.
\vspace{.1in}

\noindent
{\bf Proposition.} {\em 
For arbitrary finite set $Z$ on the unit circle,
there exists
an entire function of the form $(\ref{1})$ whose indicator diagram
coincides with the convex hull of $Z\cup-Z$.} 
\vspace{.1in}

{\em Proof.} Let $E$ be the set of all entire functions of the form
(\ref{1}). We consider the following operators on $E$:
$$R_\theta[f](z):=f(ze^{-i\theta}),$$
$$C[f](z):=\frac{1}{2}\left(f(z)+f(-z)\right),$$
and
$$S[f](z):=\frac{1}{2}\left(f(z)-f(-z)\right).$$
Now we define an operator $E\times E\to E$ by the formula
$$Q_{\theta_1,\theta_2}[f_1,f_2]=(C\circ R_{\theta_1})
[f_1]+(S\circ R_{\theta_2})[
f_2].$$
It can be easily shown that if $f\in E$ is a function with indicator
diagram $[0,1]$, then $(C\circ R_\theta)[f]$ and $(S\circ R_\theta)[f]$
have indicator diagram $[-e^{i\theta},e^{i\theta}]$.
Hence the indicator diagram of $f_1=Q_{\theta_1,\theta_2}[f,f]$
is the convex hull of $$\{ e^{i\theta_1},-e^{i\theta_1},e^{i\theta_2},
-e^{i\theta_2}\}.$$
This proves the Proposition for the sets $Z$ of two points.
Then we consider $f_2=Q_{0,\theta_3}[f_1,f]$ and so on.
\hfill$\Box$
\vspace{.1in}

Now we consider functions of the form (\ref{1}) with
$\arg a_n=2\pi n^2\alpha,\;\alpha\in\R$.
\vspace{.1in}

\noindent
{\bf Theorem 3.} {\em Let $f$ be of the form $(\ref{1})$ with
$a_n=\exp(2\pi in^2\alpha),$
where $\alpha$ is irrational.
Then $f$ has completely regular growth in the sense of Levin--Pfluger,
and $h_f\equiv 1.$}
\vspace{.1in}

We recall the main facts of the Levin--Pfluger theory in the modern
language \cite{A}. Fix a positive number $\rho$.
Let $u$ be a subharmonic function in the plane
satisfying 
$$u(z)\leq O(r^\rho),\quad r\to\infty.$$
Then the family of subharmonic functions 
$$A_tu(z)=t^{-\rho}u(tz),\quad t>1,$$
is bounded from above on every compact subset of the plane.
Such families of subharmonic functions are pre-compact in the
topology $D'$ of Schwartz's distributions \cite[Theorem 4.1.9]{H},
so from every sequence $A_{t_k}u,\, t_k\to\infty$ one can select a
convergent subsequence.
An entire function $f$ or order $\rho$, normal type is
called of {\em completely regular growth} if the limit
\begin{equation}
\label{A}
u=\lim_{t\to\infty}A_t\log|f|
\end{equation}
exists. It is easy to see that this limit is a  fixed point for all operators
$A_t$, so it has the form
$$u(re^{i\theta})=r^\rho h(\theta),$$
and $h$ is the indicator of $f$.
Operators $A_t$ also act on measures in the plane by the formula
$$A_t\mu(E)=t^{-\rho}\mu(tE)\quad\mbox{for}\quad E\subset\C.$$
The Riesz measure $\mu_f$ of $\log|f|$ is the counting measure of zeros
of $f$, and one of the results of Levin--Pfluger can be stated as follows:
The existence of the limit (\ref{A}) implies the existence
of the limit 
$$\mu=\lim_{t\to\infty}A_t\mu_f.$$
This limit $\mu$ is also fixed by all operators $A_t$,
so 
$$d\mu=r^{\rho-1}drd\nu(\theta),$$
where $\nu$ is a measure on the unit circle which is called the
{\em angular density} of zeros. This measure $\nu$ is related to
the indicator by the formula
$$d\nu=(h^{\prime\prime}+\rho^2h)d\theta,$$
in the sense of distributions.

Thus, as a corollary from Theorem~3, we obtain that the angular density
of zeros 
of $f$ is a constant multiple of the Lebesgue measure. 

Completely regular growth with indicator $1$ and order $\rho=1$ implies that
\begin{equation}
\label{7}
\log|f(re^{i\theta})|=r+o(r)\quad\mbox{as}\quad r\to\infty,
\end{equation}
uniformly with respect to $\theta$, when $re^{i\theta}$ does not belong
to an exceptional set. According to Azarin, \cite{A}, for every $\eta>0$,
this exceptional set can be covered by
discs centered at $w_k$ and of radii $r_k$ such that
\begin{equation}
\label{exp}
\sum_{|w_k|\leq r}r^\eta_k=o(r^\eta),\quad r\to\infty.
\end{equation}
This improves the original condition with $\eta=1$ given in \cite{L}.
The properties (\ref{p1}) and (\ref{p2}) of zeros of $f$, stated in
the beginning of the paper, follow from (\ref{7}) by
theorems II.2 and III.4 in \cite{L}, see also \cite{Sodin}.

The exceptional set (\ref{exp})
is larger than the exceptional set in the work of Nassif.
The exceptional set in Theorem~3 could be improved to
a set of exponentially small circles if one knew that the zeros
of $f$ are well separated. This seems to be an 
interesting unsolved problem about the function (\ref{11}).
In particular, {\em can $f$ of the form $(\ref{11})$
have a multiple zero?}
For $\alpha=\sqrt{2}$, Nassif proved that all but finitely many
zeros are simple and well separated. 

\vspace{.1in}

That the indicator of $f$ in Theorem~3 is constant was proved by 
Valiron \cite[p. 412]{Valiron}\footnote{Valiron obtained
the equation which is equivalent to our (\ref{12}) below,
\cite[Eq. (11)]{Valiron} but he did not fully explore its
consequences. Later in the same paper, on p. 421, Valiron proves
that $f$ is of completely regular growth
only under an additional Diophantine condition
on $\alpha$.}.  This also follows from the
result of Cooper \cite{Co},
who proved that the corresponding function $G$
has the unit circle as its natural boundary, see also
\cite[p.\,76, Footnote]{C}
where a short proof of Cooper's theorem is given.
However, as we noticed above,
constancy of the indicator by itself only implies
(\ref{carlson}); 
it is the statement about completely regular growth
that permits to conclude that
the zeros have positive density.
\vspace{.1in}

{\em Proof of Theorem 3.} 
By differentiating the power series
it is easy to obtain
\begin{equation}
\label{funk}
f'(z)=e^{2\pi i\alpha}f(ze^{i\beta}),\quad\mbox{where}\quad
\beta=4\pi\alpha.
\end{equation}
(This is the ``pantograph equation'' (\ref{fae}) with $q=e^{2\pi i\alpha}$.)
The assumption that $|a_n|=1$ implies the following
behavior of $M(r,f)$
\begin{equation}
\label{max}
\log M(r,f)=r+o(r).
\end{equation}
This is proved by the standard argument relating the growth
of $M(r,f)$ with the moduli of the coefficients, see, for example
\cite[Ch. I, \S2]{L}. 
The bounds 
$0\leq r-\log M(r,f)\leq (1/4+o(1))\log r$ can be obtained
as follows. The upper bound $M(r,f)\leq e^r$ is trivial, and for the lower
bound, use Cauchy's inequality $M(r,f)\geq r^n/n!$, and maximize the right
hand side with respect to $n$.
In particular, the order $\rho=1$.

It follows from (\ref{max}) that the family of subharmonic functions
$$\{ u_t= A_t\log|f|: 0<t<\infty\}$$
is uniformly bounded from above on compact subsets of $\C$.
Moreover, $u_t(0)=0$.
So every sequence $\sigma=(t_k)\to\infty$ contains a subsequence $\sigma'$
such that the limit
\begin{equation}
u=\lim_{t\in\sigma',t\to\infty}u_t
\label{conv}
\end{equation} 
exists in $D'$, the space of Schwartz's distributions in the plane.
The set of all possible limits $u$ for all sequences $\sigma$ is
called the limit set of $f$ and denoted by $\Fr[f]$.
It consists of subharmonic functions in the plane satisfying $u(0)=0$.
Equation, (\ref{max}) implies that
\begin{equation}
\label{reg}
\max_{|z|\leq r}u(z)=r,\quad 0\leq r<\infty.
\end{equation}
If $u=\lim t_k^{-1}\log|f(t_kz)|$, and $v=\lim t_k^{-1}\log|f'(t_kz)|$
with the same sequence $t_k\to\infty$, then 
\begin{equation}
\label{leq}v\leq u.\end{equation}
Indeed, by Cauchy's inequality, for every $\epsilon>0$ and
$|z|>1/\epsilon$, we have
$$\log|f'(z)|\leq\max_{|\zeta|\leq\epsilon|z|}\log|f(z+\zeta)|.$$
This implies that for every $\epsilon>0$,
$$v(z)\leq\max_{|\zeta|\leq\epsilon}u(z+\zeta).$$
Now the upper semi-continuity of subharmonic functions
shows that the right hand side of the last equation tends
to $u(z)$ as $\epsilon\to 0$, which proves (\ref{leq}).

The functional equation
(\ref{funk}) and (\ref{leq}) imply that
$u(ze^{i\beta})\leq u(z)$, and this 
gives 
\begin{equation}
\label{12}
u(ze^{i\beta})\equiv u(z).
\end{equation}
As $\beta$ is irrational, $u(z)$ is independent of $\arg z$, and
taking (\ref{reg}) into account we conclude
that the limit set $\Fr[f]$ consists of the single function $u(z)=|z|$.
This means that $f$ is of completely
regular growth with constant indicator.

\nopagebreak
\hfill$\Box$
\vspace{.1in}

Now we show that there exist irrational $\alpha$ such that
the corresponding functions $f_\alpha$ in (\ref{11}) do not
have property (\ref{ratio}).
\vspace{.1in}

\noindent
{\bf Theorem 4.} {\em There is a residual set $E$ on the unit circle,
such that for a function $f_\alpha$ as in $(\ref{11})$ with $\alpha\in E$,
we have
\begin{equation}\label{res}
\limsup_{r\to\infty}M(r,f)/m_2(r,f)=\infty.
\end{equation}
}

We recall that a set is called residual if it is an intersection
of countably many dense open sets. By Baire's Category Theorem,
residual sets on $[0,1]$ have the power of
a continuum and thus contain irrational points.
\vspace{.1in}

{\em Proof of Theorem 4.} Consider the sets
$$E_{m,n}=\{\alpha:M(r,f_\alpha)/m_2(r,f_\alpha)\leq m\quad\mbox{for}\quad
r\geq n\},$$
where $m$ and $n$ are positive integers. Evidently, all these
sets are closed. Let $E=[0,1]\backslash\cup_{m,n}E_{m,n}$. 
Then for $\alpha\in E$ we have (\ref{res}), and $E$ is a countable
intersection of open sets.
It remains to show that $E$ is dense. We will show that $E$ contains all
rational numbers. 
Indeed, for rational $\alpha$, $f_\alpha$ is a finite trigonometric
sum:
\begin{equation}
\label{cucu}
f_\alpha=\sum c_ke^{b_kz},
\end{equation}
where $b_k$ are roots of unity. This representation 
immediately follows from the functional equation (\ref{funk}):
iterating this functional equation
finitely many times, we obtain a linear differential
equation whose solutions are trigonometric sums (\ref{cucu}). 
It is clear that any finite trigonometric sum $g$ satisfies
$$M(r,g)/m_2(r,g)\to\infty\quad\mbox{as}\quad r\to\infty.$$
This proves that $E$ is dense and thus residual. 
\hfill$\Box$
\vspace{.1in}

Now we extend Theorem~3 to the case that $|a_n|$ is non-constant.
For this we need 
\vspace{.1in}

{\bf Hadamard's Multiplication Theorem} \cite{B}. {\em Let $f=\sum_{n=0}^\infty c_nz^n$ be an entire function,
and $H=\sum_{n=0}^\infty b_nz^n$ a function analytic in $\bC\backslash\{1\}$.
Then the function
$$(f\star H)(z)=\sum_{n=0}^\infty a_nc_nz^n$$
has the integral representation
$$(f\star H)(z)
=\frac{1}{2\pi i}\int_C f(\zeta)H\left(\frac{z}{\zeta}\right)
\frac{d\zeta}{\zeta},$$
where $C$ is any closed contour going once in positive direction
around the point $1$.}
\vspace{.1in}

This operation
$f\star H$ is called the Hadamard composition of power series.
{}From the integral representation we immediately obtain
\begin{equation}
\label{est}
|(f\star H)(z)| \leq K \max_{|\zeta-z|\leq r|z|}|f(\zeta)|,\quad
\mbox{where}\quad K=\max_{|\zeta-1|=r/(1-r)}|H(\zeta)|.
\end{equation}

\noindent
{\bf Theorem 5.} {\em Let $h$ be
an entire function of minimal exponential type.
Let 
\begin{equation}
\label{ass1}
c_n=h(0)h(1)\ldots h(n),\quad n\geq 0,
\end{equation}
and
assume in addition that
\begin{equation}\label{ass2}
-\log|c_n|=\frac{1}{\rho}n\log n-cn+o(n),\quad n\to\infty,
\end{equation}
with some real constant $c$.
Then the entire function
$$f(z)=\sum_{n=0}^\infty c_n e^{2\pi i n^2\alpha} z^n,$$
with irrational $\alpha$, has order $\rho$, normal type and
completely regular growth
with constant indicator.}
\vspace{.1in}

The condition that $c_{n}/c_{n-1}$ is interpolated by an entire function
of minimal exponential type is not so rigid as it may seem.
In this connection, we recall a theorem of Keldysh (see,
for example, \cite{G}) that
every function $h_1$ analytic in the sector $|\arg z|<\pi-\epsilon$
and satisfying $\log|h_1(z)|=O(|z|^\lambda),\; z\to\infty$ there,
with $\lambda=\pi/(\pi+\epsilon)<1$,
can be approximated by an entire function $h$ of normal type,
order $\lambda$ so that
$$|h(z)-h_1(z)|=O(e^{-|z|^\lambda}),\quad z\to\infty,
\quad |\arg z|<\pi-2\epsilon.$$
For example, one can take 
$$h_1(z)=z^{-1/\rho},\quad \rho>0,$$
and apply Keldysh's theorem, to
obtain a function $f$ of normal type, order $\rho$, satisfying all
conditions of Theorem~5.
\vspace{.1in}

{\em Proof of Theorem 5.} We write:
\begin{eqnarray*}
f(z)&=&\sum_{n=0}^\infty c_n e^{2\pi i n^2\alpha} z^n\\
&=&1+\sum_{n=0}^\infty c_{n+1}e^{2\pi i (n+1)^2\alpha}z^{n+1}\\
&=&1+ze^{2\pi i\alpha}\sum_{n=0}^\infty c_{n+1}e^{2\pi in^2\alpha}
\left(ze^{4\pi in\alpha}\right)^n\\
&=&1+ze^{2\pi i\alpha}(f\star H)(ze^{i\beta}),\quad\mbox{where}\quad
\beta=4\pi \alpha,
\end{eqnarray*}
and
$$H(z)=\sum_{n=0}^\infty \frac{c_{n+1}}{c_n} z^n.$$
By the assumption (\ref{ass1}) of the theorem,
and P\'olya's theorem above,
$H$ is holomorphic
in $\bC\backslash\{1\}$, so the estimate (\ref{est}) holds.
Assumption (\ref{ass2}) implies that
$$M(r,f)=\sigma r^\rho+o(r^\rho),\quad r\to\infty,$$
where $\sigma=e^{c\rho}/(e\rho),$ see \cite[Ch. I, \S2]{L}.
So from every sequence one can select a subsequence
such that the limits
$$u(z)=\lim_{k\to\infty} t_k^{-\rho}\log|f(t_kz)|\quad\mbox{and}
\quad v(z)=\lim_{k\to\infty} t_k^{-\rho}\log|(f\star H)(t_kz)|$$
exist, and (\ref{est}) implies that $v\leq u$, by a similar argument as
(\ref{leq}) was derived. 
Now the equation
\begin{equation}
\label{funk2}
f(z)=1+ze^{2\pi i\alpha}(f\star H)(ze^{i\beta})
\end{equation}
implies 
$$u(z)\leq \max\{0,v(ze^{i\beta})\}\leq\max\{0,u(ze^{i\beta})\},$$
and as $\beta$ is irrational, we conclude that $u$ does not
depend on $\arg z$. This completes the proof.
\hfill$\Box$
\vspace{.1in}

The authors thank Alan Sokal whose questions stimulated this
research, 
and Victor Katsnelson and Vitaly Bergelson for their help with literature.

{\em Purdue University, West Lafayette IN 47907-2067

eremenko@math.purdue.edu
\vspace{.1in}

Bilkent University, Ankara, Turkey}

\end{document}